\documentclass[11pt]{amsart}
\usepackage[a4paper,margin=28mm]{geometry}
\usepackage{amsmath,amssymb,mathtools}
\usepackage[colorlinks=true,linkcolor=blue,citecolor=blue,urlcolor=blue]{hyperref}
\newtheorem{theorem}{Theorem}[section]
\newtheorem{lemma}[theorem]{Lemma}
\newtheorem{proposition}[theorem]{Proposition}
\newtheorem{corollary}[theorem]{Corollary}
\newtheorem{maintheorem}{Theorem}

\theoremstyle{definition}

\newtheorem{question}[theorem]{Question}
\theoremstyle{remark}
\newtheorem{remark}[theorem]{Remark}
\numberwithin{equation}{section}
\newcommand{\Cu}{\operatorname{Cu}}
\newcommand{\Her}{\operatorname{Her}}

\newcommand{\K}{\mathcal K}
\newcommand{\F}{\mathbb F}
\newcommand{\id}{\operatorname{id}}
\newcommand{\ev}{\operatorname{ev}}
\newcommand{\N}{\mathbb N}
\newcommand{\Z}{\mathbb Z}

\newcommand{\rr}{\mathbb R}
\newcommand{\diag}{\operatorname{diag}}

\newcommand{\CE}{\mathcal E}
\newcommand{\CB}{\mathcal B}

\title[Infinite projections and pure infiniteness]{Infinite projections and failure of pure infiniteness
in simple nuclear $C^*$-algebras}

\author{Xin Ma}
	\address{X. Ma:  Institute for Advanced Study in Mathematics, Harbin Institute of Technology, Harbin, China, 150001}
\email{xma17@hit.edu.cn}

\author{Yuxuan Zheng}
\address{Y. Zheng: School of Mathematical Sciences and LPMC, Nankai University,
Tianjin, 300071, PR China.}
\email{yuxuanzheng@mail.nankai.edu.cn}

\date{\today}
\subjclass[2020]{46L05, 46L35, 46L80}
\begin{document}
\begin{abstract}
We construct a separable, unital, simple, nuclear $C^*$-algebra in the
UCT class such that every nonzero projection in its stabilization is
infinite, while it contains a nonzero projectionless hereditary
subalgebra. The construction gives a negative answer to the question
whether infiniteness of all nonzero projections in a unital simple
nuclear $C^*$-algebra implies pure infiniteness. 
\end{abstract}
\maketitle
\tableofcontents

\section{Introduction}\label{sec:introduction}

The distinction between finite and infinite projections is a basic
feature of the structure theory of operator algebras. For simple
$C^*$-algebras, however, it does not yield a type decomposition
analogous to that of von Neumann factors. R\o rdam constructed a
simple $C^*$-algebra containing both a nonzero finite projection and
an infinite projection, and obtained such an example that is
separable, nuclear and in the UCT class \cite{Rordam}. His construction
showed that a simple infinite $C^*$-algebra need not be purely
infinite, even under these additional hypotheses.

Pure infiniteness concerns the hereditary subalgebras of an algebra.
A simple $C^*$-algebra is purely infinite precisely when every nonzero
hereditary subalgebra contains an infinite projection
\cite[Proposition~6.11.5]{Blackadar}. This raises the question of
whether failure of pure infiniteness must always be witnessed by
a nonzero finite projection. Elliott posed this question in the
form recorded by R\o rdam in \cite[Question~7.7]{Rordam}. It also
appears as Problem LXXIV in the survey of Schafhauser, Tikuisis and
White \cite{STW}.

\begin{question}
Suppose that $A$ is a unital simple $C^*$-algebra and that every
nonzero projection in $A$ is infinite. Must $A$ be purely infinite?
Does the implication hold if $A$ is nuclear?
\end{question}

The implication holds whenever $A$ has property \textup{(SP)}, meaning
that every nonzero hereditary subalgebra contains a nonzero projection.
Indeed, under the projection hypothesis in the question, such a
projection is infinite. If $p$ belongs to a hereditary subalgebra
$H\subseteq A$, a partial isometry witnessing its infiniteness belongs
to $pAp\subseteq H$, so $p$ is infinite in $H$ itself. The hereditary
characterization therefore gives pure infiniteness. In particular,
the implication holds for algebras of real rank zero, which have
property \textup{(SP)}. The issue is whether infiniteness of all
nonzero projections alone rules out a nonzero projectionless
hereditary subalgebra.

Another affirmative case is supplied by $\mathcal Z$-stability, where
$\mathcal Z$ denotes the Jiang--Su algebra. Gong, Jiang and Su proved
that a unital simple $\mathcal Z$-stable $C^*$-algebra is either stably
finite or purely infinite \cite[Theorem~3]{GJS}. An infinite unit
excludes the stably finite alternative. Thus any counterexample to
the question must fail both property \textup{(SP)} and
$\mathcal Z$-stability.

Our main result provides such a counterexample within the separable
nuclear UCT class.

\begin{maintheorem}\label{main}
There exists a separable, unital, simple, nuclear $C^*$-algebra $A$
satisfying the UCT such that every nonzero projection in
$A\otimes\K$ is infinite and $A$ contains a nonzero projectionless
hereditary subalgebra. In particular, $A$ is not purely infinite.
\end{maintheorem}

The theorem answers R\o rdam's existence question affirmatively and
Problem LXXIV of \cite{STW} negatively. Since an infinite projection
in a simple $C^*$-algebra is properly infinite, all nonzero projections
in $A\otimes\K$ are in fact properly infinite. Nevertheless, they do
not occur in every nonzero hereditary subalgebra. The theorem thus
separates the infiniteness of the projections that an algebra contains
from the existence of projections in its hereditary subalgebras.

The construction starts with Suzuki's work on actions on simple
$C^*$-algebras containing finite and infinite projections \cite{Suzuki}.
We extract a unital simple nuclear algebra $B$, an action
$\alpha:\F(a,b)\curvearrowright B$, and a nonzero projection $q\in B$
such that every finite direct sum of distinct translates
$\alpha_g(q)$ is finite, although $q\oplus q$ is properly infinite.
A unital Cuntz pair $v_1,v_2\in B$ defines an injective endomorphism
\[
 \delta(x)=v_1\alpha_a(x)v_1^*+v_2\alpha_b(x)v_2^*,
 \qquad x\in B.
\]
The desired $C^*$-algebra $A$ comes from an inductive limit of matrix algebras over the mapping torus
of $\delta$ and $B$ defined in \eqref{mapping torus}, with connecting maps given by diagonal sums of forward
shifts, which is also appeared
in Kishimoto and Kumjian's construction of simple stably
projectionless $C^*$-algebras \cite[Section~2]{KK}..

\subsection*{Outline of the paper} Section~\ref{sec:preliminaries} records the comparison and
inductive-limit facts used in the proofs.
Section~\ref{sec:input} obtains the coefficient algebra from Suzuki's
construction and establishes the mapping torus as the building block.
Section~\ref{sec:limit} constructs the simple limit and proves
infiniteness of all its nonzero stabilized projections.
Section~\ref{sec:hereditary} constructs the projectionless hereditary
subalgebra and completes the proof of Theorem~\ref{main}.

\section{Preliminaries}\label{sec:preliminaries}

Throughout the paper, $\N=\{0,1,2,\ldots\}$ and $\K$ denotes the
compact operators on a separable infinite-dimensional Hilbert space.
All ideals are norm-closed and two-sided. We regard $M_k(D)$ as a
corner of $D\otimes\K$. For $a\in D_+$, we write
$\Her_D(a)=\overline{aDa}$ for the hereditary subalgebra generated
by $a$. An element is full if the closed ideal it generates is the
whole algebra. An algebra is projectionless if it has no nonzero
projections; it is stably projectionless if its stabilization is
projectionless.

For projections $p,q$, the notation $p\sim q$ means
Murray--von Neumann equivalence, and $p\precsim q$ means stable
Murray--von Neumann subequivalence. We write $p^{\oplus r}$ for the direct sum
of $r$ copies of $p$. A projection is infinite if it is equivalent
to a proper subprojection of itself, and properly infinite if
$p\oplus p\precsim p$.

For positive elements, $a\precsim b$ denotes Cuntz subequivalence:
there is a sequence $(x_j)$ in the stabilization such that
$x_j^*bx_j\to a$ in norm. These conventions agree on projections.
The Cuntz semigroup $\Cu(D)$ consists of Cuntz classes of positive
elements in $D\otimes\K$, with addition given by orthogonal direct
sum.

\begin{lemma}\label{fullmultiple}
Let $D$ be a unital $C^*$-algebra and let $p$ be a full projection
in a matrix algebra over $D$. Then $1_D\precsim p^{\oplus r}$ for some
integer $r\ge1$. If $1_D$ is infinite, then $p^{\oplus r}$ is infinite for
some $r$.
\end{lemma}
\begin{proof}
We work in the matrix algebra $A=M_k(D)$ that contains $p$.
Since $p$ is full, the linear span of $ApA$ is dense in $A$.
Choose $a_1,\ldots,a_r,b_1,\ldots,b_r\in A$ such that
$c=\sum_i a_ipb_i$ satisfies $\|c-1\|<1/2$. Put
$X=(a_1p,\ldots,a_rp)$ and $Y=(pb_1,\ldots,pb_r)^T$.
Then $c=XY$, so
\[
 0<cc^*\le\|Y\|^2XX^*=\|Y\|^2\sum_i a_ipa_i^*.
\]
As $c$ is invertible, the element $d=XX^*$ is invertible. The column
$V=X^*d^{-1/2}$ satisfies $V^*V=1$. This implies that $Q=VV^*$ is a projection. Define $P=\diag(p,\ldots,p)$, which satisfies $PV=V$. Consequently, one has $PQ=PVV^*=VV^*=Q$ and similarly $QP=Q$. This implies that $P-Q$ is a projection and thus
$VV^*\le P$. Therefore, one has $1_D\precsim 1_A\precsim P$.

On the other hand, suppose $1_D$ is infinite. Then $P=p^{\oplus r}$ contains an infinite subprojection $p_0$. Choose $v$ with $v^*v=p_0$ and $vv^*<p_0$. Then
$v+(P-p_0)$ witnesses infiniteness of $P$. 
\end{proof}

\begin{remark}\label{rem: nonunital fullmulitple}
Let $D$ be a simple $C^*$-algebra. 
\begin{enumerate}
    \item Let $p\in D$ be a full projection. Then for any projection $e\in D$, one still has $e\precsim p^{\oplus r}$ for some $r>0$. To see this, fullness provides $a_i,b_i\in D$ such that 
   \[\|e-e(\sum_ia_ipb_i)e\|\leq \|e-\sum_ia_ipb_i\|<1.\]
   Then working in the unital hereditary subalgebra $eDe$, the same argument
    of Lemma \ref{fullmultiple} shows the result. 
    \item Let $p\in D$ be an infinite projection. Then $p$ is necessarily properly infinite by applying
\cite[Proposition~6.11.3]{Blackadar} in its unital corner. 
\end{enumerate}
   
\end{remark}

\begin{lemma}\label{projectioncomparison}
Let $a\in D_+$ and let $p$ be a projection in $D\otimes\K$.
If $p\precsim a$, then $\Her_D(a)$ contains a projection
Murray--von Neumann equivalent to $p$. More generally, if
$a\in M_r(D)_+$, the equivalent projection belongs to
$\Her_{M_r(D)}(a)$.
\end{lemma}
\begin{proof}
Choose $x\in D\otimes\K$ so that
$\|p-px^*axp\|<1$. The element $d=px^*axp$ is then invertible
in the hereditary $C^*$-subalgebra $p(D\otimes \K)p$. Set $z=a^{1/2}xp\,d^{-1/2}$. We have
\[
 z^*z=p,\qquad
 zz^*=a^{1/2}xp\,d^{-1}px^*a^{1/2}\in\Her(a).
\]
Thus $zz^*$ is the required projection. Now suppose $a\in M_r(D)_+$. Then by the construction, $zz^*$ is in $M_r(D)$ as well.
\end{proof}

The following result is well-known and we omit the proof.

\begin{lemma}\label{projectionlimits}
Let $D=\overline{\bigcup_nD_n}$ be the closure of an increasing
sequence of $C^*$-subalgebras.
\begin{enumerate}
\item Every projection in $D$ is equivalent to a projection in some $D_n$.
\item If every $D_n$ is projectionless, then $D$ is projectionless.
\end{enumerate}
In particular, these hold for $D\otimes \K=\lim\limits_{\longrightarrow}M_n(D)$.
\end{lemma}

\begin{lemma}\label{simplecriterion}
Let $D=\varinjlim(D_n,\theta_n)$ be an inductive limit with injective
connecting maps. Suppose that every nonzero positive element of $D_n$ has full image in $D_m$ for some $m>n$. Then $D$ is simple.
\end{lemma}
\begin{proof}
Identify $D_n$ with their images in $D$, and let $I\lhd D$
be a nonzero ideal. Suppose $D_n\cap I=0$ holds for every $n\in \N$. Then the quotient map $\pi:D\to D/I$ is
injective on every $D_n$, it would be isometric on their dense
union, hence on $D$. This would contradict $I\ne0$. Thus
$I\cap D_n\ne0$ for some $n$. Choose $0\ne a\in(I\cap D_n)_+$.
For every $k\ge n$, regard $a$ as a nonzero positive element of
$D_k$ and choose $m\ge k$ such that $a$ is full in $D_m$. Note that $D_m\cap I$ is a closed ideal in $D_m$ containing $a$. Then the fullness of $a$ implies that $D_m\cap I=D_m$, and thus $D_m\subseteq I$. Furthermore, one has $D_k\subseteq I$ and thus $I=D$. 
\end{proof}

\section{Construction of the mapping torus}\label{sec:input}
In this section, we construct the mapping torus $E$ as the building block of the final $C^*$-algebra $A$.

\begin{lemma}\label{lem: key1}
There exist a separable, unital, simple, nuclear $C^*$-algebra $B$
satisfying the UCT, an action $\alpha:\F(a,b)\curvearrowright B$, a nonzero projection $q\in B$,
and isometries $v_1,v_2\in B$ with
\[
 v_i^*v_j=\delta_{ij}1,\qquad v_1v_1^*+v_2v_2^*=1,
\]
such that for every finite subset $F\subset\F(a,b)$, the projection
\[
 q_F=\bigoplus_{g\in F}\alpha_g(q)
\]
is finite in $B\otimes\K$.
The unit of $B$ is properly infinite and $[1_B]_{K_0(B)}=0$.
The projection $q\oplus q$ is properly infinite.
\end{lemma}

\begin{proof}
We  follow the construction in the proof of 
\cite[Theorem 1]{Suzuki} for $\Gamma=\F(a,b)$, the free group with $2$ generators.  Define
\[
 C=c_0(\Gamma)\otimes C(Z)\otimes\K,
 \qquad Z=(S^2)^{\N}.
\]
Let $\mathcal E=\varinjlim(\mathcal M (C), \varphi)$ be the inductive limit $C^*$-algebra defined in \cite[Equation (14)]{Suzuki} and denote its canonical maps by $\mu_j: \mathcal M(C)\to \mathcal E$. Define
$D=C^*(\mu_j(C):j\in\Z)$, and writes $R=D\rtimes_\sigma\Z$, which is simple by \cite[Claim 13]{Suzuki}. Denote by $\eta$ the action of $\Gamma$ on $R$ described in \cite[Paragraphs before Claim 10]{Suzuki}.
Let $p^B\in M_2(C(S^2))$ be the Bott projection. Write $\pi_j: Z\to S^2$ the $j$-th coordinate projection. Denote by $p_0=\pi^*_0(p^B)\in M_2(C(Z))\subset C(Z)\otimes\K$ and write
$e_{11}\in\K$ be a rank-one projection. With $\delta_e$ denoting
the characteristic function of the group identity, set
\[
 q_0=\mu_0(\delta_e\otimes p_0),\qquad
 r_0=\mu_0(\delta_e\otimes1_{C(Z)}\otimes e_{11}).
\]
Both are projections in $D\subset R$. Let
\[
 \Phi:R\longrightarrow\mathcal M(D\otimes\K(\ell^2(\Z)))
\]
be the faithful regular representation used in \cite[Claim~14]{Suzuki}, and let
$f\in\K(\ell^2(\Z))$ be the rank-one projection in 
\cite[equation (31)--(32)]{Suzuki}, which imply that, in the multiplier algebra, for 
$Q=\bigoplus_{g\in\Gamma}\Phi(\eta_g(q_0))$, one has
\[
 r_0\otimes f\not\precsim Q,
 \text{ and } r_0\otimes f\precsim\Phi(r_0).
\]
Note that by explicit calculation, $\eta_g(q_0)=\mu_0(\delta_{g^{-1}}\otimes p_0)$.

For a finite set $F\subset\Gamma$, suppose
$P_F=\bigoplus_{g\in F}\eta_g(q_0)$ is infinite. Since $R\otimes\K$ is simple, it would be full and properly
infinite.  It follows from Remark \ref{rem: nonunital fullmulitple} that 
\begin{equation*}\label{eq1}
    r_0\precsim P^{\oplus N}_F\precsim P_F
\end{equation*}
for some $N$ by fullness  and proper infiniteness of $P_F$.
 Applying $\Phi$ would give
\[
 r_0\otimes f\precsim\Phi(r_0)
 \precsim\bigoplus_{g\in F}\Phi(\eta_g(q_0))\precsim Q,
\]
which is a contradiction. Thus all $P_F$ are finite.


Choose $p$ as in the final paragraph of the proof of
\cite[Theorem~1]{Suzuki}: $p$ is a nonzero properly infinite projection
with $[p]_{K_0(R)}=0$, and there is a nonzero properly infinite
projection orthogonal to $p$. The existence of such projections uses
\cite{Cuntz}. Suzuki's argument provides multiplier unitaries carrying
$\eta_a(p)$ and $\eta_b(p)$ to $p$. We obtain an action $\alpha$ of
$\Gamma$ on $B=pRp$ by cocycle perturbing $\eta$ on the free generators
$a,b$.
Explicitly, the action has the form
\[
 \alpha_g(x)=u_g\eta_g(x)u_g^*,\qquad
 u_g\eta_g(p)u_g^*=p,
\]
where $u_g\in\mathcal M(R)$. More precisely, choose the unitaries on $a,b$ as in Suzuki's
argument, set $u_{s^{-1}}=\eta_{s^{-1}}(u_s^*)$ for $s=a,b$, and
extend along reduced words by $u_{gh}=u_g\eta_g(u_h)$ by the freeness of $\Gamma$.

Because unital $B$ is a full corner of $R$, then $B$
is simple, separable and
nuclear because $R$ is. In addition, $B$ satisfies the UCT because it is strongly Morita equivalent to $R$. On the other hand,
$p$, as the unit of $B$, is sill properly infinite in $B$. Moreover, because $B$ is the full corner, the full-corner inclusion induces an isomorphism $K_0(B)\xrightarrow{\cong}K_0(R)$ and therefore, $[p]_{K_0(B)}=0$ holds as well.

Using the proper infiniteness of $p$, Remark \ref{rem: nonunital fullmulitple} implies that $q_0\precsim p$. Then there is a projection $q\leq p$ that is equivalent to $q_0$. Choose a $v\in R$ such that $v^*v=q_0$ and $vv^*=q\in B$. Define 
$t_g=u_g\eta_g(v)\in R$, which satisfies
\[
 t_g^*t_g=\eta_g(q_0),\qquad t_gt_g^*=\alpha_g(q).
\]
Thus, one has $q_F=\bigoplus_{g\in F}\alpha_g(q)$ is equivalent to $\bigoplus_{g\in F}\eta_g(q_0)=P_F$ in $R\otimes\K$. Because $P_F$ is shown finite in $R\otimes \K$, so is $q_F$ in $R\otimes \K\simeq B\otimes \K$.

We now show $q\oplus q$ is properly infinite. First, it follows from \cite[Claim 11]{Suzuki} that $r_0$ is
properly infinite in $R$. In addition, \cite[Equation (26)]{Suzuki} implies $r_0\precsim q_0\oplus q_0$. Thus $q_0\oplus q_0$ is infinite,
and hence properly infinite by simplicity. Its equivalent copy
$q\oplus q$ has the same property in $M_2(B)$.

Finally, we construct desirable isometries $v_1, v_2$. We work in $M_2(B)$ and note that 
\[[\begin{pmatrix}1_B&0\\0&1_B\end{pmatrix}]=[\begin{pmatrix}1_B&0\\0&0\end{pmatrix}]=0\]
in $K_0(B)$. Both projections are infinite in the simple unital
algebra $M_2(B)$, so \cite[Corollary~6.11.9]{Blackadar} gives
$1_B\oplus1_B\sim1_B$. Let the equivalence is implemented by $W\in M_2(B)$ with
\[
W^*W=\begin{pmatrix}1&0\\0&1\end{pmatrix},\qquad
 WW^*=\begin{pmatrix}1&0\\0&0\end{pmatrix}.
\]
Then observe $\diag(1_B,0)W=WW^*W=W(1_B\oplus 1_B)=W$. This implies that $W$ necessarily has the form
\[W=\begin{pmatrix}v_1&v_2\\0&0\end{pmatrix}.\]
Then direct calculation implies that 
\[
 v_i^*v_j=\delta_{ij}1,\qquad v_1v_1^*+v_2v_2^*=1
\]
as we want.
\end{proof}

Now let $B$, $\Gamma=\F(a,b)$, $q, v_1, v_2\in B$ satisfy Lemma \ref{lem: key1}.
Define
\begin{equation}\label{delta}
\delta(x)=v_1\alpha_a(x)v_1^*+v_2\alpha_b(x)v_2^*,\qquad x\in B.
\end{equation}
Note that the relation $v_i^*v_j=\delta_{ij}1$ and $v_1v_1^*+v_2v_2^*=1$ imply that 
$\delta$ is a unital
$*$-homomorphism. In addition, 
since
$v_1^*\delta(x)v_1=\alpha_a(x)$, one has $\delta(x)=0$ implying $\alpha_a(x)=0$, which means $x=0$. Thus $\delta$ is isometric.

Define the mapping torus
\begin{equation}\label{mapping torus}
   E=\{f\in C([0,1],B):f(1)=\delta(f(0))\} 
\end{equation}
which is separable and unital. The extension
\[
 0\longrightarrow C_0((0,1),B)\longrightarrow E
 \xrightarrow{\operatorname{ev}_0} B\longrightarrow0
\]
shows that $E$ is nuclear and satisfies the UCT: both the ideal
and the quotient have these properties, and the nuclear UCT class
is closed under extensions.

\begin{lemma}\label{shifts}
There is a point-norm continuous semigroup $(\rho_s)_{s\ge0}$
of unital injective endomorphisms of $E$ with $\rho_0=\id_E$.
For every $s\ge0$, $\rho_s$ is homotopic to $\id_E$. For
$f\in M_k(E)$ and $t\in[0,1]$, it satisfies
\[
 \|\rho_s^{(k)}(f)(t)\|=\|f(\{t+s\})\|,
\]
where $\{u\}=u-\lfloor u\rfloor$ and endpoint norms are understood
periodically. Then for any projection $P\in M_k(E)$, one has
\[ \rho_s^{(k)}(P)\sim P.\]
\end{lemma}
\begin{proof}
For $f\in E$, extend $f$ to $[0,\infty)$ by
\[
 \widetilde f(t+k)=\delta^k(f(t))
 \quad(0\le t<1,\ k\in\mathbb N).
\]
The function $\widetilde f$ is continuous on every open interval $(n, n+1)$. At any integer $n\geq 1$, the limit from left is 
\[\delta^{n-1}(f(1))=\delta^{n-1}(\delta(f(0)))=\delta^n(f(0)),\]
which implies that $\widetilde f$ is continuous on $[0, \infty)$. Set
\[
 \rho_s(f)(t)=\widetilde f(t+s),\qquad s\ge0,\quad0\le t\le1.
\]
For every $u\ge0$ and integer $k\ge0$, observe that
\[
 \widetilde f(u+k)=\delta^k(\widetilde f(u)),
\]
which entails $\rho_s(f)(1)=\delta(\rho_s(f)(0))$ and thus $\rho_s(f)\in E$. Moreover, a direct calculation shows that each $\rho_s$ is a unital $*$-homomorphism. Let $t+s=u+n$ for $u=\{s+t\}\in [0,1)$ and $n\in \N$. Then 
\[\rho_s(f)(t)=\delta^n(f(u)).\]
Recall $\delta$ is isometric, so is $\delta^n$. This means
\[\|\rho_s(f)(t)\|=\|f(u)\|\]
and therefore
\[\|\rho_s(f)\|=\sup_{t\in [0,1]}\|f(\{s+t\})\|=\|f\|.\]
This shows that $\rho_s$ is injective. The same argument yield
\[\|\rho^{(k)}_s(f)(t)\|=\|f(\{s+t\})\|.\]

We now claim $(\rho_s)_{s\geq 0}$ form a semigroup. Indeed, for $u=t+n$ with $t\in [0, 1)$ and $n\in \N$, one has
\[\widetilde{\rho_s(f)}(u)=\delta^n(\rho_s(f)(t))=\delta^n(\widetilde f(t+s))=\widetilde f(u+s).\]
Therefore, one has
\[\rho_r(\rho_s(f))(t)=\widetilde{\rho_s(f)}(t+r)=\widetilde f(t+r+s)=\rho_{r+s}(f)(t).\]
This shows
\[
 \rho_s\rho_t=\rho_{s+t},\qquad\rho_0=\id_E.
\]

We finally shows $\rho_s$ is homotopic to $\id_E$. Let $s_0, s$ be given. Note that
\[\|\rho_{s_0}(f)-\rho_s(f)\|=\sup_{t\in [0,1]}\|\widetilde f(s_0+t)-\widetilde f(s+t)\|.\]
Therefore, if $s\to s_0$, then $\|\rho_s(f)-\rho_{s_0}(f)\|\to 0$. This implies $s\mapsto\rho_s(f)$ is norm continuous on every bounded
interval. Thus $r\mapsto\rho_{rs}$, for $0\le r\le1$, is a point-norm
homotopy from the identity to $\rho_s$. Similar argument shows $r\mapsto \rho^{(k)}_{rs}$ for $0\le r\le 1$ is a point-norm homotopy from the identity to $\rho^{(k)}_s$.

Finally let $P\in M_k(E)$ be a projection. Then for $\lambda\in [0,1]$, define $P_\lambda=\rho^{(k)}_{\lambda s}(P)$, which form a continuous path from $P$ to $\rho^{(k)}_s(P)$. Then choose a finite increasing sequence $0=\lambda_0<\cdots<\lambda_l=1$ such that $P_{\lambda_0}=P$, $P_{\lambda_l}=\rho_s^{(k)}(P)$, and 
$\|P_{\lambda_i}-P_{\lambda_{i+1}}\|<1$. Then $P_{\lambda_i}$ is unitary equivalent to $P_{\lambda_{i+1}}$. This implies that 
\begin{equation}\label{projshift}
 \rho_s^{(k)}(P)\sim P.
\end{equation}
\end{proof}

\begin{lemma}\label{full}
An element $f\in M_k(E)_+$ is full if $f(t)\neq 0$ for any $t\in [0,1]$.
Consequently every nonzero projection in $M_k(E)$ is full.
\end{lemma}
\begin{proof}
Denote by $\mathcal B=M_k(B)$ and $\Delta=\delta^{(k)}$ for simplicity. Then 
\[\mathcal E:=M_k(E)=\{h\in C([0,1],\mathcal B): h(1)=\Delta(h(0))\}.\]
Let $I$ be the ideal generated by $f$. Define $J=C_0((0,1), \mathcal B)$, which can be identified as $\{h\in \mathcal E: h(1)=h(0)=0\}$, which is a closed two-sided ideal in $\mathcal E$. Define $K=I\cap J$.

Let $t_0\in (0,1)$ be arbitrary and choose a scalar function $\chi\in C_c((0,1))$ with $\chi(t_0)=1$. Note that $\chi 1_{\mathcal B}\in \mathcal E$ and therefore $\chi f=\chi 1_{\mathcal B}\cdot f\in I$ because $I$ is an ideal. Moreover, by definition, $\chi f\in J$ and therefore $\chi f\in K$ satisfying 
\begin{equation}\label{eq2}
   (\chi f)(t_0)=f(t_0)\neq 0. 
\end{equation}
For the evaluation map $\ev_{t_0}: J\to \CB$ at $t_0$, it is not hard to see that $\ev_{t_0}(K)$ is a closed $2$-sided ideal in $\CB$. Since $\CB$ is simple, it follows from \eqref{eq2} that $\ev_{t_0}(K)=\CB$.

Now, let $h\in C_c((0,1), \CB)$. Write $L=\operatorname{supp}(h)$ for simplicity. Let $\epsilon>0$. For any $t\in L$, because $\ev_t(K)=\CB$, there exists a $k_t\in K$ such that $k_t(t)=h(t)$. By continuity, there is an open neighborhood $U_t\subset (0,1)$ of $t$ on which one has
\begin{equation}\label{eq3}
    \|k_t(u)-h(u)\|<\epsilon.
\end{equation}
Then the compactness of $L$ yields a finite subcover of $\{U_t: t\in L\}$, say, $U_1,\dots, U_m$, with corresponding sections $k_1,\dots, k_m\in K$. Choose $\varphi_i\in C_c(U_i)_+$ such that
$\sigma=\sum_{i=1}^m\varphi_i\le1$ everywhere and $\sigma=1$ on $L$, and define
\[k=\sum_{i=1}^m\varphi_ik_i.\]
This belongs to $K$ as each $\varphi_ik_i\in K$. Since $\sigma h=h$,
\[
 h-k=\sum_{i=1}^m\varphi_i(h-k_i),
\]
and \eqref{eq3} gives $\|h-k\|<\epsilon$. Since $\epsilon>0$ is arbitrary, and $K$ is closed, one has $h\in K$, which implies $C_c((0,1), \CB)\subset K$ and thus $J=K$. This further shows $J\subset I$.

Note that the evaluation map $\ev_0: \CE\to \CB$ is surjective with the kernel $J$. On the other hand, one also has $\ev_0(I)=\CB$. Thus, for an arbitrary $x\in \CE$, there exists $y\in I$ such that $x(0)=y(0)$. This implies that $x-y\in \ker(\ev_0)=J\subset I$. So, one has $x\in I$ and thus $\CE=I$.

For a projection $P$, the continuous function $t\mapsto\|P(t)\|$
takes values in $\{0,1\}$. On the connected interval it is constant.
Thus a nonzero projection is nonzero at every fiber.
\end{proof}

\begin{lemma}\label{infiniteunit}
The unit of $E$ is infinite.
\end{lemma}
\begin{proof}
Put $w=v_2\alpha_b(v_2)$. We have $w^*w=1$ and $v_1^*w=0$. The other required
orthogonality follows from
\[
 \begin{aligned}
 \delta(v_1)^*w
 &=\bigl(v_1\alpha_a(v_1^*)v_1^*
        +v_2\alpha_b(v_1^*)v_2^*\bigr)v_2\alpha_b(v_2)\\
 &=v_2\alpha_b(v_1^*v_2)=0.
 \end{aligned}
\]
Define
\[
 z(t)=
 \begin{cases}
 \cos(\pi t)v_1+\sin(\pi t)w,&0\le t\le\tfrac12,\\
 \cos(\pi(t-\tfrac12))w+
 \sin(\pi(t-\tfrac12))\delta(v_1),&\tfrac12\le t\le1.
 \end{cases}
\]
Orthogonality gives $z(t)^*z(t)=1$ throughout. Moreover
$z(1)=\delta(z(0))$, which implies $z\in E$. Its range projection is proper because $z(0)z(0)^*=v_1v_1^*<1$. Hence $1_E$ is infinite.
\end{proof}

\begin{corollary}\label{stageamplification}
For every nonzero projection $P$ in a matrix algebra over $E$,
there is an integer $N\ge1$ such that $P^{\oplus N}$ is infinite.
\end{corollary}
\begin{proof}
Lemma~\ref{full} makes $P$ full. By Lemma~\ref{fullmultiple},
some multiple $P^{\oplus N}$ of $P$ dominates the infinite projection $1_E$
from Lemma~\ref{infiniteunit}. Therefore $P^{\oplus N}$ is infinite.
\end{proof}

\section{The inductive limit}\label{sec:limit}
In this section, we construct the $C^*$-algebra $A$ the a inductive limit of matrices over $E$ and prove that all non-zero projections in $A\otimes\K$ is infinite.

Let $S_0=\{0\}$, $N_0=1$, and choose integers $m_n\ge2$ tending to
infinity. Recursively choose an integer
\[
 L_n>\operatorname{diam}(S_{n-1})+3,
\]
and define
\[
 \tau_{n,j}=jL_n+\frac{j}{m_n}\quad(0\le j<m_n),\qquad
 S_n=S_{n-1}+\{\tau_{n,j}:0\le j<m_n\},\qquad
 N_n=N_{n-1}m_n.
\]
Define $A_n=M_{N_n}(E)$ and unital injective maps
\[
 \varphi_n:A_{n-1}\longrightarrow A_n,
 \qquad f\longmapsto
 \operatorname{diag}_{0\le j<m_n}
       \big(\rho_{\tau_{n,j}}^{(N_{n-1})}(f)\big).
\]
Write $A=\varinjlim(A_n,\varphi_n)$. For $m>n$ let
$\varphi_{m,n}=\varphi_m\circ\cdots\circ\varphi_{n+1}$,
and let $\varphi_{n,n}=\id_{A_n}$. We identify each stage with
its image in $A$ whenever convenient. This is possible because
the connecting maps are injective.

\begin{lemma}\label{separation}
For every $n\ge0$, the set $S_n$ has cardinality $N_n$ and its
distinct elements are at distance greater than $2$. The map
$\varphi_{n,0}:E\to A_n$ is, up to a  permutation of diagonal
blocks, of the form $f\mapsto\diag_{s\in S_n}\rho_s(f)$.
\end{lemma}
\begin{proof}
The assertions are immediate for $n=0$. Suppose they hold for $n-1$.
If $s,s'\in S_{n-1}$ and two new sums use the same index $j$,
their distance is $|s-s'|>2$. If their indices $j,j'$ differ, the
reverse triangle inequality gives
\[
 |s+\tau_{n,j}-s'-\tau_{n,j'}|
 \ge L_n-\operatorname{diam}(S_{n-1})-1>2.
\]
Thus all $N_{n-1}m_n=N_n$ sums are distinct and satisfy the gap
bound.

For the final statement,  First for $N_0=1$ and $S_0=\{0\}$, note that 
\[\varphi_{1,0}(f)=\varphi_1(f)=\diag_{0\leq j_1< m_1}\rho_{\tau_{1,j_1}}(f)\]
with $S_1=\{\tau_{1,j_1}: 0\leq j_1<m_1\}$. Then at the second stage,  since $(\rho_s)_{s\geq 0}$ is a semigroup by Lemma \ref{shifts}, one has
\begin{align*}
    \varphi_{2,0}(f)&=\diag_{0\leq j_2<m_2}\left(\rho_{\tau_{2,j_2}}[\diag_{0\leq j_1<m_1}\rho_{\tau_{1,j_1}}(f)]\right)\\
    &=\diag_{0\leq j_2<m_2}\diag_{0\leq j_1<m_1}\rho_{\tau_{2,j_2}+\tau_{1,j_1}}(f)
\end{align*}
with $S_2=\{\tau_{1,j_1}+\tau_{2, j_2}: 0\leq j_1<m_1, 0\leq j_2<m_2\}$. Then by induction, one actually has
\[\varphi_{n,0}(f)=\diag_{0\leq j_n<m_n}\dots\diag_{0\leq j_1<m_1}\rho_{\tau_{1,j_1}+\dots+\tau_{n, j_n}}(f)\]
and $S_n=\{\tau_{1,j_1}+\dots+\tau_{n, j_n}: 0\leq j_k<m_k, k=1,\dots,n\}$. Note that we have shown that any two elements in $S_n$ are different to each other. This implies that any $s\in S_n$ appeared only once among the diagonal entries in $\varphi_{n, 0}(f)$, which implies
$\varphi_{n, 0}(f)=\diag_{s\in S_n}\rho_s(f)$
up to a permutation of diagonal blocks.
\end{proof}

Consequently, for every fixed $t\in[0,1]$, the integers
\begin{equation}\label{distinct}
 \lfloor t+s\rfloor\quad(s\in S_n)
\end{equation}
are distinct. Composition of the connecting maps sends an original
section $f\in E$ to the block diagonal section
\begin{equation}\label{paths}
 \operatorname{diag}_{s\in S_n}\rho_s(f),
\end{equation}
up to an irrelevant permutation of matrix blocks.

\begin{proposition}\label{simple}
The limit $A$ is unital, separable, simple and nuclear, and it
satisfies the UCT.
\end{proposition}
\begin{proof}
All stages are unital and separable, and all connecting maps are
unital. Nuclearity and the UCT pass from $E$ to its matrix algebras
and to the inductive limit. It is left to show the simplicity, for which we apply

Fix $0\ne f\in(A_n)_+=M_{N_n}(E)_+$. By Lemma~\ref{simplecriterion}, it suffices to show $\varphi_{k,n}(f)$ is full in $A_k$.
Define a continuous function $g: t\mapsto\|f(t)\|$ on $[0,1]$. Since $\delta^{(N_n)}$ is isometric, one has
\[g(1)=\|\delta^{(N_n)}(f(0))\|=\|f(0)\|=g(0).\]
This induce a continuous function, still denoted by $g$, on the circle $\mathbb T=\rr/\Z$. Because $f\neq 0$, there exists a $t_0$ such that $g(t_0)>0$. By continuity, choose an open
arc $U$ on which $g$ is positive. Normalizing the length of the circle $\mathbb T$ by $1$ and we denote by $l$ the length of the arc $U$.

For $k>n$, define
$f_{k-1}=\varphi_{k-1,n}(f)$. Recall that each $\tau_{r, 0}=0$ and thus connecting map $\varphi_r$ is of the form
\[\varphi_r(x)=\diag(x, \rho^{(N_{r-1})}_{\tau_{r,1}}(x),\dots \rho^{(N_{r-1})}_{\tau_{r, m_r-1}}(x)).\]
So $f_{k-1}$ contains the original section
$f$ as a diagonal block. Hence
$\|f_{k-1}(t)\|\ge\|f(t)\|$ for every $t\in [0,1]$.

Choose $k>n$ with $1/m_k$ smaller than the length $l$ of $U$.
For every $t\in[0,1]$, the set
\[\{t+j/m_k: j=0,\dots, m_k-1\}\]  
has to meet $U$, say $t+j/m_k\in U$. Then it follows from Lemma \ref{shifts} that 
\[
 \|\rho^{(N_{k-1})}_{\tau_{k,j}}(f_{k-1})(t)\|
 =\|f_{k-1}(\{t+j/m_k\})\|
 \ge\|f(\{t+j/m_k\})\|>0.
\]
Thus $\varphi_{k,n}(f)$ is nonzero at every fiber and is full by
Lemma~\ref{full}. Thus, $A$ is simple.
\end{proof}

\begin{proposition}\label{allproj}
Every nonzero projection in $A\otimes\K$ is infinite.
\end{proposition}
\begin{proof}
Let $0\ne p\in A\otimes\K$ be a projection. By
Lemma~\ref{projectionlimits}, it is equivalent to the image of a
nonzero projection $P\in M_k(A_n)=M_{kN_n}(E)$ for some $k,n$. Write $d=kN_n$ for simplicity.
Then, it follows from Corollary~\ref{stageamplification} that there exists an integer $N$ such that
$P^{\oplus N}$ is infinite.

For $m>n$, set $R=R_{m,n}:=\prod_{j=n+1}^m m_j$. Define
\[T_{m,n}=\{\tau_{n+1, j_{n+1}}+\dots+\tau_{m,j_m}: 0\leq j_l< m_l, l=n+1,\dots, m\}.\]
The sums are distinct by Lemma~\ref{separation}, after prefixing the
zero choices at stages $1,\ldots,n$. Thus $|T_{m,n}|=R_{m,n}$.
Enumerate $T_{m,n}=\{s_1,\ldots,s_R\}$. Then, up to a permutation, one has
\[\varphi_{m,n}^{(k)}(P)=\diag(\rho^{(d)}_{s_1}(P),\dots, \rho^{(d)}_{s_R}(P))\in M_{Rd}(E)=M_k(A_m).\]
Now Lemma \ref{shifts} shows that 
\[\rho_{s_r}^{(d)}(P)\sim P\]
inside $M_d(E)$ for each $1\leq r\leq R$. Let $w_r$ witness that $w_r^*w_r=P$ and $w_rw_r^*=\rho_{s_r}^{(d)}(P)$. Define $W=\diag(w_1,\dots,w_R)$, which satisfies $W^*W=P^{\oplus R}$ and \[WW^*=\diag(\rho^{(d)}_{s_1}(P),\dots,\rho^{(d)}_{s_R}(P))\]
and thus 
\[
 \varphi_{m,n}^{(k)}(P)\sim P^{\oplus R}
 \quad\text{in }M_k(A_m).
\]
Choose $m$ with $R=R_{m,n}\ge N$. Then $P^{\oplus R}$ contains
infinite subprojection $P^{\oplus N}$ and thus it is infinite itself. Thus, $\varphi^{(k)}_{m,n}(P)$ is infinite. Denote by $\iota_i^{(k)}: M_k(A_i)\to A\otimes \K$ the canonical embedding. Then $p$ is equivalent to $\iota^{(k)}_n(P)=\iota_m^{(k)}(\varphi^{(k)}_{m,n}(P))$, which is infinite. Thus $p$ is infinite.
\end{proof}

\section{A projectionless hereditary subalgebra in $A$}\label{sec:hereditary}
In this section, we construct a projectionless hereditary subalgebra in $A$ and prove the main theorem.
\begin{lemma}\label{fixedfinite}
Let $B$ be a simple unital algebra with infinite
unit, let $\delta:B\to B$ be a unital endomorphism, and suppose
that a projection $q\in B$ has the property that
$\bigoplus_{d\in F}\delta^d(q)$ is finite for every finite set
$F\subset\N$. Let $D\subset\N$ be nonempty and finite, and write
\[
 Q_D=\bigoplus_{d\in D}\delta^d(q).
\]
There is no nonzero projection $p\le Q_D$ satisfying
$p\sim\delta^{(|D|)}(p)$, where $\delta^{(|D|)}$ is the matrix amplification of $\delta$.
\end{lemma}
\begin{proof}
Suppose such a $p$ exists. Applying $\delta$ to the equivalence $p\sim\delta^{(|D|)}(p)$
repeatedly gives
\[
 p\sim(\delta^k)^{(|D|)}(p)\qquad(k\ge0).
\]
Choose an integer $M>\max D-\min D$. For each $r\ge1$, the sets
$D+iM$ for $0\le i<r$, are pairwise disjoint. Therefore
\[
 \bigoplus_{i=0}^{r-1}(\delta^{iM})^{(|D|)}(Q_D)
   =\bigoplus_{\substack{0\le i<r\\d\in D}}\delta^{iM+d}(q)
\]
is finite by the hypothesis on distinct iterates, since all the
integers $iM+d$ are distinct. Its subprojection
$\bigoplus_{i=0}^{r-1}(\delta^{iM})^{(|D|)}(p)$ is equivalent to
$p^{\oplus r}$. Thus every finite multiple of $p$ is finite.

On the other hand, $p$ is nonzero and hence full in $M_{|D|}(B)$. Consequently $1_B\precsim p^{\oplus r'}$ for
some $r'$ and thus $p^{\oplus r'}$ is infinite by Lemma \ref{fullmultiple}. This is a
contradiction.
\end{proof}

For the rest of the section, let $B$ be the $C^*$-algebra defined in Lemma \ref{lem: key1}. In addition,
let $\alpha: \Gamma=\F(a,b)\curvearrowright B$ be the action obtained in Lemma \ref{lem: key1}. In addition, let $q, v_1, v_2\in B$ be elements in Lemma \ref{lem: key1}. We still denote by $\delta$ the isometric unital $*$-homomorphism defined in \eqref{delta}.
Choose $0\ne h\in C_0((0,1))_+$ and set
\[
 c(t)=h(t)q\in E_+.
\]
Write $c_n$ for its image in $A_n$, and $c_\infty$ for its nonzero
image in $A$. Let $H_n=\Her_{A_n}(c_n)$.

\begin{lemma}\label{wordlemma}

For every $k\ge0$,
\begin{equation}\label{word}
 \delta^k(q)\sim
 \bigoplus_{w\in\{a,b\}^k}\alpha_w(q).
\end{equation}
Consequently, for every finite set $D\subset\N$ of distinct
integers, $\bigoplus_{d\in D}\delta^d(q)$ is finite.
\end{lemma}
\begin{proof}

Write $v_a=v_1$ and $v_b=v_2$.  Set $V_\varnothing=1$ for the
empty word. For
$w=i_1\cdots i_k$ with $i_j\in\{a,b\}$, set
\[
 V_w=v_{i_1}\alpha_{i_1}(v_{i_2})\cdots
          \alpha_{i_1\cdots i_{k-1}}(v_{i_k}).
\]
Induction on $k$, using the Cuntz relations for $v_1, v_2$ obtained in Lemma \ref{lem: key1} and multiplicativity
of the automorphisms, one actually verifies $V_w^*V_{w'}=\delta_{w,w'}1$ for
words of the same length, and
\[
 \sum_{|w|=k}V_wV_w^*=1,\qquad
 \delta^k(q)=\sum_{|w|=k}V_w\alpha_w(q)V_w^*.
\]
Define $T_k=(V_w\alpha_w(q))_{|w|=k}\in M_{1,2^k}(B)$. The $(w,w')$-entry of $T^*_kT_k$ is $\alpha_w(q)V_w^*V_{w'}\alpha_{w'}(q)$, which is zero if $w\neq w'$ and is $\alpha_w(q)$ if $w=w'$. Denote by $\mathcal W_k=\{a,b\}^k$ the set of positive words of length $k$. This implies that 
$T_k^*T_k=\diag_{w\in \mathcal W_k}\alpha_w(q)$ and 
\[T_kT_k^*=\sum_{w\in \mathcal W_k}V_w\alpha_w(q)V_w^*=\delta^k(q).\]
This has established \eqref{word}. For the final statement, let $D\subset \N$ be finite.  Observe that 
\[\bigoplus_{d\in D}\delta^d(q)\sim\bigoplus_{d\in D}\bigoplus_{w\in \mathcal W_d}\alpha_w(q),\]
which is finite by Lemma \ref{lem: key1}.
\end{proof}

Fix $n$ and $t\in[0,1]$. By \eqref{paths},
$c_n=\diag_{s\in S_n}\rho_s(c)$ and thus
\[c_n(t)=\operatorname{diag}_{s\in S_n}\left(h(\{t+s\})\delta^{\lfloor t+s\rfloor}(q)\right).\]
It is supported by the projection
\[
 Q_n(t)=\operatorname{diag}_{s\in S_n}
              \delta^{\lfloor t+s\rfloor}(q)
\]
in the sense that $Q_n(t)c_n(t)=c_n(t)Q_n(t)=c_n(t)$.
We remark that no continuity of $t\mapsto Q_n(t)$
is asserted here. By \eqref{distinct} and \eqref{word}, $Q_n(t)$
is equivalent to a finite direct sum of \emph{distinct} translates
of $q$. It is therefore finite by Lemma~\ref{lem: key1}.

\begin{proposition}\label{herfinite}
Every $H_n$ is projectionless. Consequently
$H=\Her_A(c_\infty)$ is a nonzero projectionless hereditary subalgebra.
\end{proposition}
\begin{proof}
Let $P\in H_n$ be a projection. If $P\ne0$, then each $P(t)$ is a projection and $t\mapsto \|P(t)\|$ is continuous and thus constant. Therefore, $p=P(0)\ne0$. The path
$t\mapsto P(t)$ of projections in $M_{N_n}(B)$ has endpoints
\[
 p=P(0),\qquad P(1)=\delta^{(N_n)}(p).
\]
Hence $p\sim\delta^{(N_n)}(p)$.

Moreover, since $P\in H_n$, there exists a sequence $a_j$ such that $c_na_jc_n\to P$ and thus $c_n(0)a_j(0)c_n(0)\to P(0)=p$. Recall that $c_n(0)$ is supported on $Q_n(0)$ and therefore 
$p\le Q_n(0)=\bigoplus_{s\in S_n}\delta^{\lfloor s\rfloor}(q)$.
The integers $\lfloor s\rfloor$ are distinct by \eqref{distinct}. Denote by $D_n=\{\lfloor s\rfloor: s\in S_n\}$ such that $|D_n|=|S_n|=N_n$. Note that $p$ satisfies $p\leq Q_n=Q_{D_n}=\bigoplus_{d\in D_n}\delta^{d}(q)$ and $p\sim \delta^{(|D_n|)}(p)$, which is a contradiction to Lemma~\ref{fixedfinite} Thus $H_n$
has no nonzero projections.

Finally, identifying the stages through the injective maps, density gives
\[
 H=\overline{c_\infty A c_\infty}
   =\overline{\bigcup_n\overline{c_n A_n c_n}}
   =\overline{\bigcup_n H_n}.
\]
Projectionlessness then passes to this increasing union by Lemma \ref{projectionlimits}(2).
\end{proof}

Now we are ready to prove Theorem \ref{main}.

\begin{proof}[Proof of Theorem~\ref{main}]
Proposition~\ref{simple} shows that $A$ is unital separable nuclear simple and satisfy the UCT. Then
Proposition~\ref{allproj} shows infiniteness of every nonzero 
projections in $A\otimes \K$. However, Proposition~\ref{herfinite} supplies a nonzero
projectionless hereditary subalgebra. Recall that a simple purely infinite
$C^*$-algebra has an infinite projection in every nonzero hereditary
subalgebra by \cite[Proposition~6.11.5]{Blackadar}. Therefore $A$ is
not purely infinite. 
\end{proof}

\section{Acknowledgement}
During  the development of this work, the authors used ChatGPT 6 Astra for searching necessary references, examining mathematical arguments, and language editing.  The authors have checked and substantially revised all mathematical arguments generated by AI and take full responsibility for the final version.

\end{document}